\documentclass[a4paper]{amsart}

\usepackage[english]{babel}
\usepackage[utf8]{inputenc}
\usepackage{amsmath,amssymb}
\usepackage[all]{xy}
\usepackage{graphicx}
\usepackage[colorinlistoftodos]{todonotes}

\newtheorem{theorem}{Theorem}
\newtheorem{corollary}[theorem]{Corollary}

\newtheorem{lemma}[theorem]{Lemma}
\newtheorem{proposition}[theorem]{Proposition}

\newtheorem{example}{Example}[section]

\newcommand{\cO}{{\mathcal{O}}}
\newcommand{\cG}{{\mathcal{G}}}
\newcommand{\cF}{{\mathcal{F}}}

\newcommand{\cB}{{\mathcal{B}}}

\newcommand{\Q}{{\mathbb{Q}}}
\newcommand{\Z}{{\mathbb{Z}}}

\title[On the computation of the Ap\'ery set]{On the computation of the Ap\'ery set of numerical monoids and affine semigroups}

\author{Guadalupe M\'arquez--Campos}
\address{Depto. de \'Algebra and IMUS, Universidad de Sevilla. E-- 41080 Sevilla, Spain.}
\email{gmarquez@us.es}
\author{Ignacio Ojeda}
\address{Depto. de Matem\'{a}ticas, Universidad de Extremadura,  E--06071 Badajoz, Spain.}
\email{ojedamc@unex.es}
\author{Jos\'e M. Tornero}
\address{Depto. de \'Algebra and IMUS, Universidad de Sevilla. E--41080 Sevilla, Spain.}
\email{tornero@us.es}
\thanks{First and third authors were partially supported by the grant FQM--218 and P12--FQM--2696 (FEDER and FSE). The second author is partially supported by the project  MTM2012-36917-C03-01, National Plan I+D+I and by Junta de Extremadura (FEDER funds)}

\date{\today}

\begin{document}

\begin{abstract}
A simple way of computing the Ap\'ery set of a numerical semigroup (or monoid) with respect to a generator, using Groebner bases, is presented, together with a generalization for affine semigroups. This computation allows us to calculate the type set and, henceforth, to check the Gorenstein condition which characterizes the symmetric numerical subgroups.
\end{abstract}

\maketitle

\section{Introduction}

The most elementary structure this paper deals with is that of numerical monoid. A numerical monoid is a very special kind of semigroup that can be thought of as a set 
$$
\langle \, a_1,...,a_k \, \rangle = \left\{ \lambda_1 a_1 + ... + \lambda_k a_k \; | \; \lambda_i \in \Z_{\geq 0} \right\}, \mbox{ with } \gcd(a_1,...,a_k)=1.
$$

This object has been thoroughly studied in the last years, when a significant number of problems concerning the description of these semigroups and some of their more interesting invariants have been tackled. Unless otherwise stated, all proofs which are not included can be found in \cite{GS,RA}.

Given a numerical monoid $S$, there are some invariants which will be of interest for us. The most relevant will be the set of gaps, noted $G(S)$, and defined by
$$
G(S) = \Z_{\geq 0} \setminus S,
$$
which is a finite set. Its cardinal will be noted $g(S)$ and its maximumm, noted $f(S)$, is called the Frobenius number of $S$.

The Ap\'ery set of $S$ \cite{Apery} with respect to an element $s \in S$ can be defined as 
$$
Ap(S,s)= \{0, w_0, . . . ,w_{s-1}\}
$$
where $w_i$ is the smallest element in $S$ congruent with $i$ modulo $s$. 

This set has some highly interesting properties, which we will review now. 

\begin{lemma}\label{lema0}
With the previous notations, for all $s \in S$, 
$$
Ap(S,s)= \{ x \in S \; \mid \; x - s  \notin S \}
$$
\end{lemma}

\begin{example} 
Let $S= \langle 7, 9, 11, 15 \rangle$. Some of the Ap\'ery sets associated to its generators are:
$$
Ap(S, 7) =\{0, \, 9,\, 11,\, 15,\, 20,\, 24,\, 26\}
$$
$$
Ap(S, 15) = \{0,\, 7,\, 9,\, 11,\, 14,\, 16,\, 18,\, 20,\, 21,\, 23,\, 25,\, 27,\, 28,\, 32,\, 34\}
$$
\end{example}

In particular, for a monoid with two generators, the Ap\'ery sets associated to its generators are fully determined.

\begin{lemma}
Let $S = \langle \, a_1, \; a_2 \, \rangle$. Then
$$
Ap(S,a_i) = \left\{0, \, a_j, \, 2a_j, ..., \, (a_i-1)a_j \right\}
$$
\end{lemma}

The importance of the Ap\'ery set can be illustrated in the following result, where the previous numerical invariants of a numerical monoid are expressed in terms of Ap\'ery sets.

\begin{proposition}
Let $S$ be a monoid, $s \in S$. Then:
$$
f(S) = \max \big\{Ap(S,s) - s \big\}, \quad \quad g(S) = \frac{1}{s} \sum_{w \in Ap(S,s)} w + \frac{s-1}{2}.
$$
\end{proposition}

This expression of $g(S)$ in terms of the Ap\'ery set is known as Selmer's formula.

Our first purpose is to develop an easy algorithm to compute the Ap\'ery set of a $S$ associated to any of its generators. As it turns out, these are the most interesting ones. We will do that with the help of Groebner bases with respect to some monomial orderings, and more specifically in terms of the characterization of elements in (respectively, out) a monoid in terms of these bases.

This characterization can be found in \cite{MCT}, but we will review it in the second section of the paper for the convenience of the reader (Theorem \ref{GC}). It will be extensively used in the sequel as it provides an easy (not necessarily in complexity terms) and versatile characterization of elements in (or not in) $S$, with a quite wide choice of monomial orderings. 

After that, in the third section, we show how to compute (and {\em see}) the Ap\'ery set associated to the generators (Theorem \ref{Apery}). We are interested not only on the computation, but also on the variety of choices we can make for the orderings. This will be the most important feature in the sequel.

The computation can be generalized to a wider class of semigroups: the affine semigroups. This will be done in the fourth section of the paper at the cost of being more restrictive with the orderings.

The computation of the Ap\'ery set will prove useful in order to study a related concept: the set of pseudo--Frobenius elements ($PF(S)$) and the type set. We will explore this relationship in the fifth section. Precisely, if we consider, for a given $S = \langle a_1,...,a_k \rangle$, with $a_{1}<\cdots<a_{k}$, its type set 
$$
T(S) = \left\{ m \in Ap(S,a_k) \; | \; m+a_i \notin Ap(S,a_k), \; \forall i=1,...,k \right\},
$$
then we have \cite{NW,Kunz} 
$$
g(S) = \frac{1}{2} \left( f(S)+1 \right) \; \Longleftrightarrow \; \# T(S) = 1.
$$

This condition ($\# T(S) =1$) always holds for monoids with two generators \cite{Sylvester}, but in general it is not so. The condition is subsequently known as the {\em Gorenstein condition}. Concerning this, we will show how we can compute the type set by calculating the Ap\'ery set $Ap(S,a_k)$ for a variety of monomial orderings (Theorem \ref{GS}).

We will end with a word on pseudo--Frobenius elements and how they are related with more complex algebraic structures which can be defined from $S$.

Finally, something should be said on previous work. The Ap\'ery set has been thoroughly studied in the literature \cite{CJZ, CZ, RR, MCH, RGGB} but its computation has not been treated so often. 

The work by P. Pis\'on--Casares \cite{Pilar} used a technique close to Groebner bases (although the orderings used there were not strictly monomial orderings) to compute Ap\'ery sets in the context of affine semigroups. Our third and fourth section are actually a twist of her arguments, in order to use monomial orderings. 

In particular, our main result in section $3$ (generalized in section $4$) is included anyway, both for the convenience of the reader and also because it provides a detailed description of the several orderings that we can consider. This will be the most important thing in section $6$ and hopefully in future applications.

\section{The Groebner correspondence}

All the general results referred to the general theory of Groebner bases that we will use can be found in standard textbooks on the subject; for instance in \cite{AL,CLOS}. 

We will work in the polynomial ring $\Q[x, y_1, ..., y_k]$. Given a polynomial 
$$
g = \sum_{\underline{\alpha} = (\alpha_0,...,\alpha_k) \in \Z^{k+1}_{\geq 0}} a_{\underline{\alpha}} x^{\alpha_0} y_1^{\alpha_1} ... y_k^{\alpha_k} \in \Q[x, y_1, ..., y_k],
$$
and a monomial ordering $\prec$, define
$$
\exp (g) = \max_{\prec} \left\{ \underline{\alpha} \; | \; a_{\underline{\alpha}} \neq 0 \right\} \in \Z^{k+1}_{\geq0}.
$$

Let then $S = \langle a_1,...,a_k \rangle \subset \Z_{\geq 0}$ be a numerical monoid. Consider the following binomial ideal associated to $S$:
$$
I_S = \langle y_1 - x^{a_1},  y_2 - x^{a_2} ,  y_3 - x^{a_3} , ... ,  y_k - x^{a_k} \rangle \subset \Q[x, y_1, ..., y_k].
$$

As $I_S$ is a binomial ideal (i.e., it can be generated by binomials), it has a lot of special properties \cite{Binomial}. In particular we are interested in the following ones:

\begin{itemize}
\item A reduced Groebner basis of $I_S$ consists of binomials.
\item The normal form of a monomial with respect to a reduced Groebner basis of $I_S$ is again a monomial.
\end{itemize}

So, consider $I_S$ and let $\cB =\{g_1,...,g_r\}$ be the reduced Groebner basis of $I_S$ with respect to $\prec$, an elimination ordering for $x$, and $N_\cB$ the normal form with respect to $\cB$.

Let us write also:
$$
q_i = \exp_\prec (g_i)\in \Z^{k+1}_{\geq0}, \quad\quad K_{q_i} = q_i + \Z^{k+1}_{\geq0}, \quad\quad E \left( I_S \right) = \bigcup_{i=1}^r K_{q_i} \subset \Z^{k+1}_{\geq0}.
$$

A preliminary elementary result is:

\begin{lemma}\label{GC0}
Let $\widetilde{\phi}$ be the ring homomorphism defined by: 
\begin{eqnarray*}
\widetilde{\phi} : \Q[x, y_1, y_2, ..., y_k] & \longrightarrow & \Q[x] \\
y_j & \longmapsto & x^{a_j}\\
x & \longmapsto & x
\end{eqnarray*}

Then $\ker \left( \widetilde{\phi} \right) = I_S$. Therefore
$$
N_\cB (g) = N_\cB (h) \; \Longleftrightarrow \; g-h \in \ker \left( \widetilde{\phi} \right) 
\; \Longleftrightarrow \; \widetilde{\phi} (g) = \widetilde{\phi} (h).
$$
\end{lemma}

The main result in \cite{MCT} goes as follows: 

\begin{theorem}\label{GC}
With the previous notations, let $\overline{E(I_S)} = \Z^{k+1} \setminus E(I_S)$. Then:
\begin{itemize}
\item The mapping
\begin{eqnarray*}
\mathcal{F}: G(S) &  \longrightarrow & \overline{E(I_S)} \setminus \{ x =0 \} \subset \Z^{k+1}_{\geq0}\\
N & \longmapsto & \exp \left(N_{\cB} \left(x^{N} \right) \right)
\end{eqnarray*}
is bijective.
\item The mapping
\begin{eqnarray*}
\mathcal{G} : S & \longrightarrow & \overline{E(I_S)} \bigcap \{ x =0 \} \subset \Z^{k+1}_{\geq0} \\
M & \longmapsto &  \exp \left( N_\cB \left(x^{M} \right) \right)
\end{eqnarray*}
is bijective.
\end{itemize}
Furthermore, for any $l \in \Z_{\geq 0}$, if 
$$
N_\cB \left(x^l\right) = \left( \sigma_0,\sigma_1,...,\sigma_k \right),
$$
then 
$$
l = \sigma_0 + \sigma_1 a_1 + ... + \sigma_k a_k.
$$
\end{theorem}

So the mapping $l \longmapsto N_\cB \left(x^l\right)$ separates the non--negative integers into (or out of) $\overline{E(I_S)} \cap \{x=0\}$, depending on whether they are in $S$ or not.

\begin{example} 
Let us consider an example of dimension $3$ (taken from \cite{MCT}). Let $S=  \langle 7,9,11 \rangle$. The Frobenius number of this numerical semigroup is:
$$
f(S) =26,
$$
and its set of gaps:
$$
G(S) = \{1,\ 2,\ 3,\ 4,\ 5,\ 6,\ 8,\ 10,\ 12,\ 13,\ 15,\ 17,\ 19,\ 24,\ 26\}.
$$

We can take the binomial ideal:
$$
I_S = \langle y_1 - x^{7},  y_2 - x^{9}, y_3 - x^{11}  \rangle \subset \Q[x, y_1, y_2, y_3]
$$
and find the Groebner basis $\cB$, using an elimination ordering for $x$. For this example, we have taken the usual lexicographic ordering $x > y_1 > y_2 > y_3$. With this particular choice we get:
\begin{eqnarray*}
\mathcal{B} &=& \big\{ \ y_{2}^{11}-y_{3}^9, \; -y_{2}^2+y_{3}y_1, \; y_{2}^9y_1-y_{3}^8, \; y_{2}^7 y_{1}^2-y_{3}^7, \; y_{2}^5 y_{1}^3-y_{3}^6, \; y_{2}^3 y_{1}^4-y_{3}^5, \\
    & &  \quad y_{1}^5 y_{2}-y_{3}^4, \; -y_{2} y_{3}^3+y_{1}^6, \; -y_{2} y_{1}^2+y_{3}^2 x, \; -y_{1}^3+y_{3} y_{2} x, \; y_{2}^3 x-y_{1}^4,  \\
    & & \quad y_{2}^2 y_{1}^2 x-y_{3}^3, \; -y_{3}^2+y_{1}^3 x, \; y_{2} x^2-y_{3}, \; y_{1} x^2-y_{2}, \; y_{3} x^3-y_{1}^2, \; -y_{1}+x^7 \ \big\}.
\end{eqnarray*}

We have to consider then, $q_{i} = \exp(g_i)$ where $g_i$ is the $i$--th polynomial in $\cB$, and take the corresponding set
$$
K_{q_{i}} = q_{i} + \Z_{\geq0} ^{k+1} \subset \Z^{k+1}_{\geq 0},
$$
in order to establish our bijections $\cF$ and $\cG$. In this case,
$$
\begin{array}{rclrclrcl}
q_1 &=& (0,0,11,0), \;\; & q_2 &=& (0,1,0,1), \; \; & q_3 &=& (0,1,9,0), \;\; \\
q_4 &=& (0,2,7,0), & q_5 &=& (0,3,5,0), & q_6 &=& (0,4,3,0), \\
q_7 &=& (0,5,1,0), & q_8 &=& (0,6,0,0), & q_9 &=& (1,0,0,2), \\
q_{10} &=& (1,0,1,1), & q_{11} &=& (1,0,3,0), & q_{12} &=& (1,2,2,0), \\
q_{13} &=& (1,3,0,0), & q_{14} &=& (2,0,1,0), & q_{15} &=& (2,1,0,0), \\
q_{16} &=& (3,0,0,1), & q_{17} &=& (7,0,0,0)
\end{array}
$$

Let us have a closer look to $\mathcal{F}$, so we are only interested in points of $\overline{E(I)}$ outside $x=0$. In order to represent the points, we will consider the subcases $x = \lambda$, with $\lambda \in \Z_{\geq 0}$. We have then:

\begin{itemize}
\item $x=1$. In this hyperplane we find several corners $q_i$, precisely
$$
q_9=(0,0,2), \; q_{10}=(0,1,1), \; q_{11}=(0,3,0), \; q_{12}= (2,2,0), \; q_{13}=(3,0,0)
$$

\noindent These points determine the elements of $\overline{E(I_S)} \subset \Z^4_{\geq 0}$, along with the point $(1,1,0,1) \in K_{q_2}$. In the following pictures we will draw square points for points in $\cup K_{q_i}$ and round points for points outside $\cup K_{q_i}$, thus associated with a unique element of $G(S)$ by means of $\cF$:

\vspace{.5cm}

\begin{center}
\includegraphics[scale=0.4]{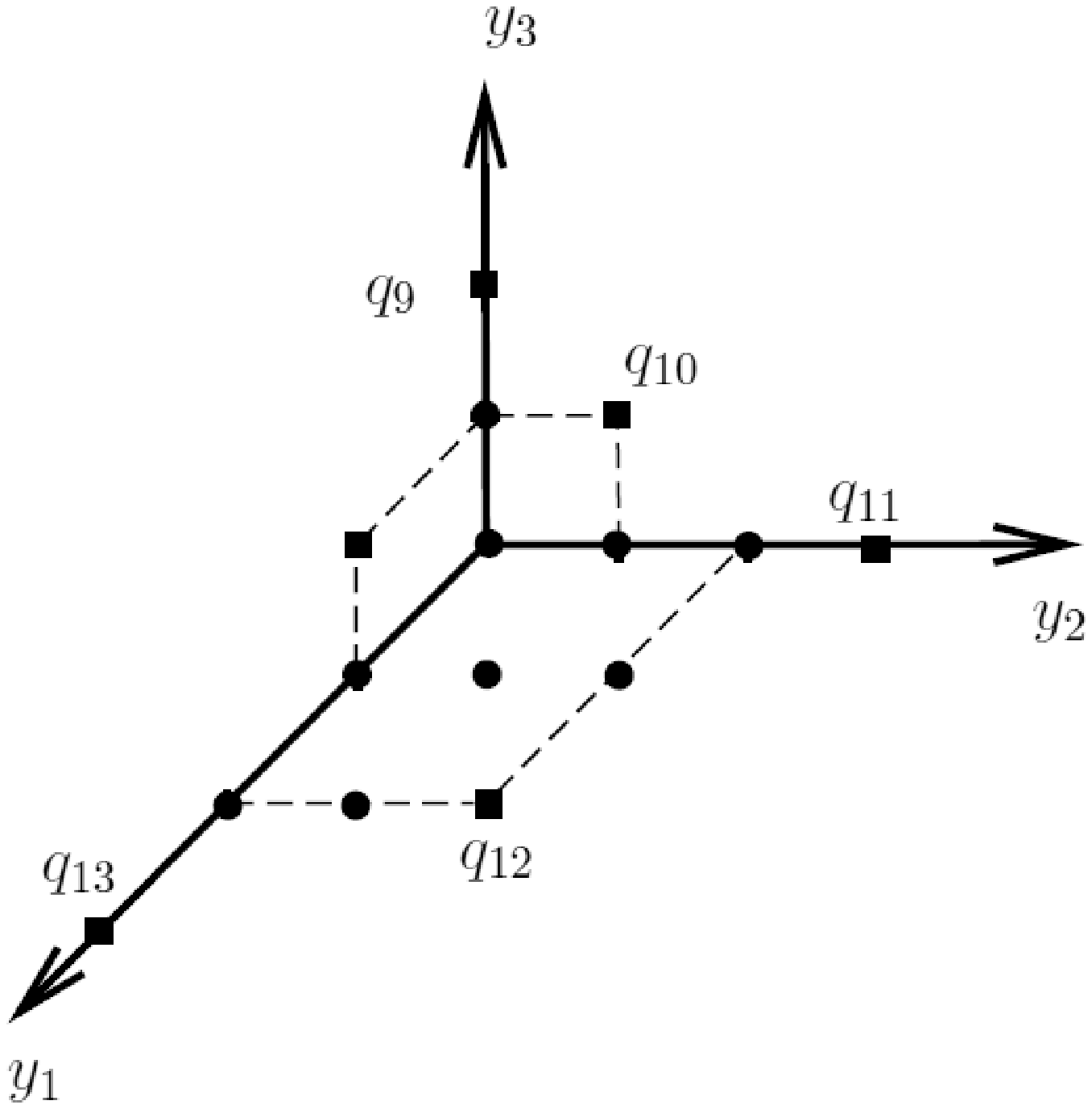}
\end{center}

\vspace{.5cm}

\item At $x=2$ these are the points which determine the set:
$$
q_{14}=(0,1,0), \; \; q_{15}=(1,0,0), \; \; (0,0,2) \in K_{q_9}
$$

\vspace{.5cm}

\begin{center}
\includegraphics[scale=0.4]{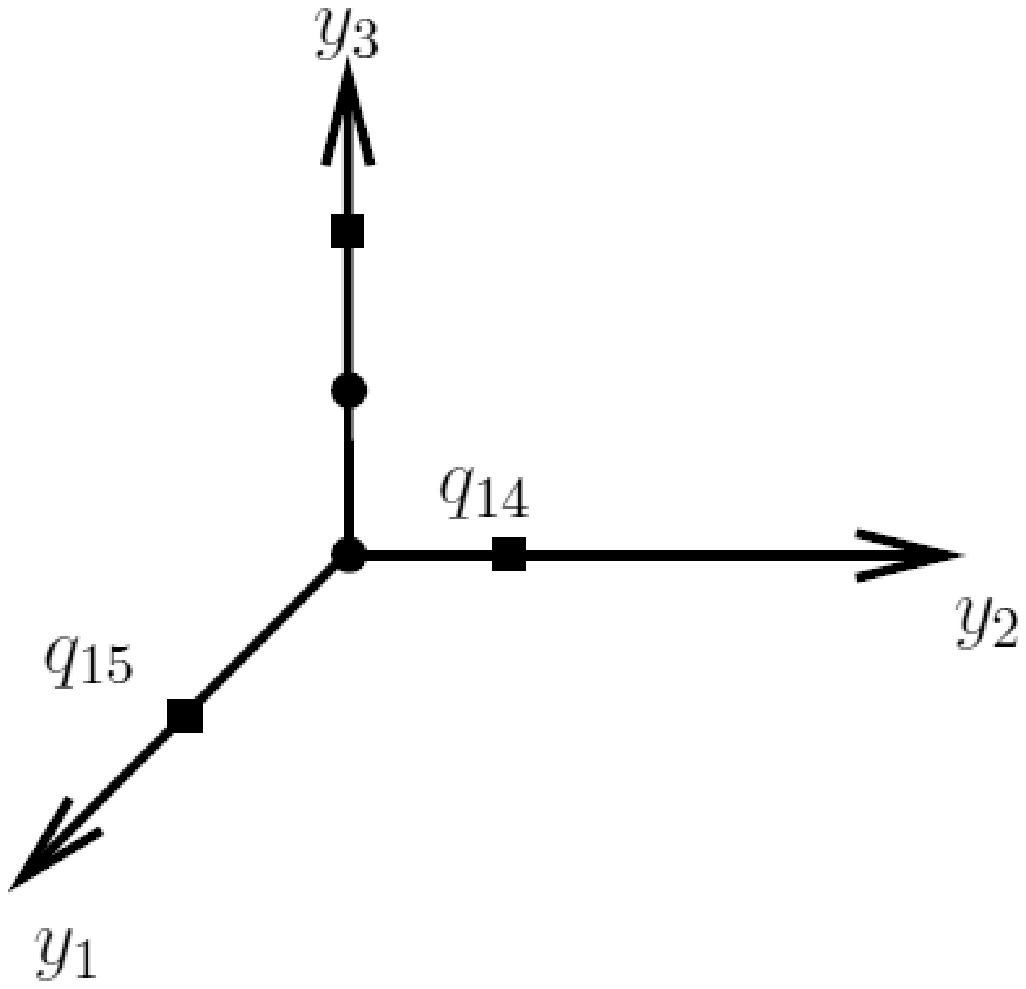}
\end{center}

\vspace{.5cm}

\item At $x=3$, we have these points in $\cup K_{q_i}$
$$
q_{16}=(0,0,1), \; \; (1,0,0) \in K_{q_{15}}, \; \; (0,1,0) \in K_{q_{14}}
$$

\vspace{.5cm}

\begin{center}
\includegraphics[scale=0.4]{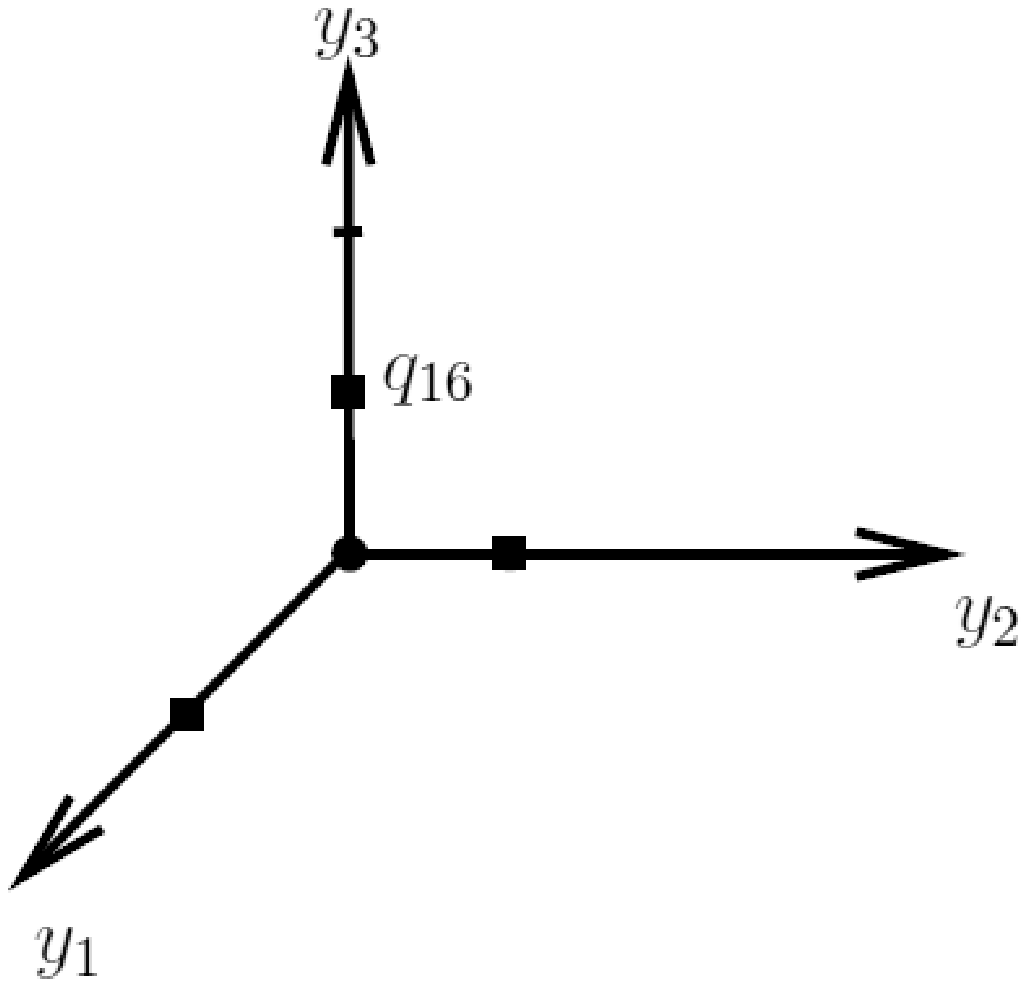}
\end{center}

\vspace{.5cm}

\item At $x=4$, $x=5$ and $x=6$, the only relevant point is the origin, as $y_i < 1 $ for $i=1,2,3$.

\item Last, in $x=7$ we have $(7,0,0,0) = q_{17}$, so this is, so to speak, the \emph{ceiling} for variable $x$.
\end{itemize}

If we compute the normal form of monomials $x^{n_i}$, where $n_i$ is the $i$--th gap, we get:
$$
\begin{array}{lcllcllcl}
N_{\mathcal{B}} (x^1) &=& x \quad & N_{\mathcal{B}} (x^2) &=& x^2  \quad & N_{\mathcal{B}} (x^3) &=& x^3 \\
N_{\mathcal{B}} (x^4) &=& x^4 & N_{\mathcal{B}} (x^5) &=& x^5 & N_{\mathcal{B}} (x^6) &=& x^6 \\
N_{\mathcal{B}} (x^8) &=& x y_1 & N_{\mathcal{B}} (x^{10}) &=& x y_2 & N_{\mathcal{B}} (x^{12}) &=& x y_3 \\
N_{\mathcal{B}} (x^{13}) &=& x^2 y_3 & N_{G} (x^{15}) &=& x y_1^2 & N_{\mathcal{B}} (x^{17}) &=& x y_1 y_2 \\
N_{\mathcal{B}} (x^{19}) &=& x y_2^2 \quad & N_{\mathcal{B}} (x^{24}) &=& x y_1^2 y_2 \quad & N_{\mathcal{B}} (x^{26}) &=& x y_1 y_2^2
\end{array}
$$
\end{example}

\section{Computation of the Ap\'ery set: the numerical case}

Consider the polynomial ring $\Q[x, y_1, ..., y_k]$, and let us define an elimination ordering for $x$, written $\prec_j$, as follows:

\begin{enumerate}
\item First, we take into account the exponent on $x$.

\item After that, we take a graded ordering in $\{ y_1, ..., y_{k} \}$ leaving aside $y_j$, and determined by the generators $a_1,...,a_k$, that is, we order by 
$$
\sum_{i=1,i\neq j }^{i=k} \alpha_i a_i,
$$ 
where $\alpha_i$ is the exponent on $y_i$.

\item Then we use any monomial ordering for all variables $y_i$, where $i=1,...,k$ and $i \neq j$.

\item Finally, we use the exponent on $y_j$.
\end{enumerate}

We will call such an ordering an Ap\'ery ordering (with respect to $a_j$).

\begin{example} 
If one wants to use the matrix notation for the ordering, an example of $\prec_j$ will be given by\footnote{Note that, if one needs a square matrix, as some computer algebra packages do, one can always add a numb variable or erase a row between the $3$rd and the $(k+1)$--th.}
$$
\left(
\begin{array}{ccccccc}
  1 &  0  & 0   & ... & \overbrace{0}^{(j+1)} & ... & 0\\
  0 & a_1 & a_2 & ... & 0 & ... & a_k\\
  0 &  1  &  0  & ... & 0 & ... & 0\\
  0 &  0  &  1  & ... & 0 & ... & 0\\
  \vdots & \vdots & \vdots & & \vdots & & \vdots \\
  0 &  0  &  0  & ... & 0 & ... & 1\\
  0 &  0  &  0  & ... & 1 & ... & 0\\
\end{array}
\right)
$$

This will be precisely the matrix of the Ap\'ery ordering we should use if we consider in (3) the lexicographic ordering $y_1<...<y_{j-1}<y_{j+1}
<...<y_k$.
\end{example}

Let $S= \langle a_1,...,a_k  \rangle$ be a numerical monoid with $a_i \neq 0$, and the ideal $I_S$ defined as in the previous section. Let $\cB_j$ be the reduced Groebner basis of $I_S$ with respect to $\prec_j$, and let us write $N_j$ the normal form with respect to this basis.

We define the following set:
$$
\Delta_{\prec_j} (S,a_j) = \left\{ \, N \in \Z_{\geq 0} \; | \; \exp \left( N_{j} \left( x^N \right) \right) \in \{x = y_j =0\}  \cap \overline{E(I_S)} \, \right\}.
$$

The important result in this section is the following:

\begin{theorem}\label{Apery}
$\Delta_{\prec_j} (S,a_j)= Ap(S, a_j)$.
\end{theorem}

\begin{proof}
Let us show first $Ap(S, a_j) \subseteq \Delta_{\prec_j}(S,a_j)$. Let $n \in Ap(S,a_j)$ other than $0$ (as the zero case is trivial), which implies $n \in S$ and $n > a_j$. There are $x_1,...,x_k \in \Z_{\geq 0}$ such that
$$
n=\sum_{i=1}^{k} a_i x_i
$$
and, being in the Ap\'ery set, we already know $n-a_j \notin S$.

We want to prove 
$$
\exp \left(N_j \left(x^n \right)\right) \in \{x=0\} \cap \{y_j = 0\}.
$$ 

But the exponent lies in $\{x=0\}$ from Theorem \ref{GC}. Let us write 
$$
\exp \left( N_j \left( x^n \right) \right)= (\gamma_1, \gamma_{2},..., \gamma_k).
$$

From the expression above
$$
\begin{array}{rcl}
n-a_j &=& a_1 \gamma_1 +  ...+ a_{j} \gamma_{j} + ...+ a_{k} \gamma_{k} - a_j \\
&=& a_1 \gamma_1 + ...+ a_j (\gamma_j -1) +...+ a_{k} \gamma_{k}.
\end{array}
$$

As  $n-a_j \in G(S)$, the above expression must have a strictly negative coefficient. As $\gamma_i \in \Z_{\geq 0}$ for $i =1,...,k$; it must be $(\gamma_j -1) \notin \Z_{\geq 0}$ and therefore $\gamma_j =0$, as we wanted to show.

Let us prove now $\Delta_{\prec_j} (S,a_j) \subseteq Ap(S,a_j)$, so take $n \in \Delta_{\prec_j} (S,a_j)$. From the definition of $\Delta_{\prec_j} (S,a_j)$ we know $n \in S$ and 
$$
\exp \left( N_j \left( x^n\right) \right) \in \{y_j =0\} \cap \{x =0\}.
$$

Hence there must exist
$$
\gamma_1,...,\gamma_{j-1},\gamma_{j+1}, ...,\gamma_{k} \in \Z_{\geq 0},
$$ 
such that
$$
N_j \left( x^n \right)= y_1^{\gamma_1}...y_{j-1}^{\gamma_{j-1}} \cdot y_{j+1}^{\gamma_{j+1}}...y_{k}^{\gamma_{k}},
$$
which, from Theorem \ref{GC} yields
$$
n = \gamma_1 a_1 + ... \gamma_{j-1}a_{j-1} + \gamma_{j+1}a_{j+1} + ... + \gamma_k a_k.
$$

Let us assume $n \notin Ap(S,a_j)$. So, either $n \notin S$ (which is impossible as $n \in \Delta_{\prec_j} (S,a_j)$) or $n > a_j$ and besides $n- a_j \in S$.

Then there are $\alpha_1,...,\alpha_k \in \Z_{\geq 0}$ verifying
$$
n-a_j = \sum_{i=1}^{k} a_i \alpha_i,
$$
that is,
$$
n=\sum_{i=1, i \neq j}^{k} a_i \alpha_i + a_j (\alpha_j +1).
$$

So we have two expressions for $n$, and
$$
\begin{array}{rcl}
n &=& a_1 \gamma_1 + ... +a_{j-1} \gamma_{j-1} + a_{j+1} \gamma_{j+1} +...+ a_{k} \gamma_{k} \\
&=& a_1 \alpha_1 + ...+a_j ( \alpha_j +1) +... + a_k \alpha_k,
\end{array}
$$
which yields
$$
\sum_{i=1, i \neq j}^{k}  a_i(\alpha_i- \gamma_i) + a_{j}(\alpha_j +1)= 0.  \quad \quad \quad \quad (*)
$$

From Lemma \ref{GC0} and the definition of $(\gamma_1,...,\gamma_k)$,
$$
y_1^{\gamma_1} \cdot ... \cdot y_{j-1}^{\gamma_{j - 1}} \cdot y_{j+1}^{\gamma_{j + 1}} \cdot...\cdot y_{k}^{\gamma_{k}}
= N_j \left( x^n \right) = N_j \left( y_1^{\alpha_1} \cdot ... \cdot y_j^{\alpha_j + 1} \cdot ... \cdot y_{k}^{\alpha_{k}} \right)
$$

So, as the normal form of a monomial must be a monomial, we know that
$$
(0, \gamma_1, ...,\gamma_{j-1},0,\gamma_{j+1},..., \gamma_{k}) \prec_j (0,\alpha_1, \alpha_2, ..., \alpha_k)
$$
which implies
$$
\sum_{i=1, i \neq j}^{k} \gamma_i a_i  \leq \sum_{i=1,i \neq j }^{k} \alpha_i a_i \; \Longrightarrow \;
\sum_{i=1, i \neq j}^{k} a_i \left( \alpha_i - \gamma_i \right) \geq 0.  
$$

This contradicts $(*)$, as $a_j (\alpha_j +1)$ must be positive.
\end{proof}

Note that our assumptions on the monomial ordering are in fact quite sparse. This will give us a lot of elbow room to work with, and we will take advantage of this fact.

\begin{corollary}
The set $\Delta_{\prec} (S,a_j)$ does not depend on the choice of $\prec$ (as long as it is an Ap\'ery ordering).
\end{corollary}

\begin{example} Let $S = \langle 7,8,9,13 \rangle$. We have
$$
I_S = \left\langle y_1 - x^{7},  y_2 - x^{8} ,  y_3 - x^{9}, y_4 - x^{13} \right\rangle \subset \Q[x, y_1, y_2, y_3,y_4],
$$

For this monoid we have
$$
G(S) = \{ \; 1,\ 2,\ 3,\ 4,\ 5,\ 6,\ 10,\ 11,\ 12,\ 19 \; \}
$$
and
$$
Ap(S, 13) = \{ \; 0,\ 7,\ 8,\ 9,\ 14,\ 15,\ 16,\ 17,\ 18,\ 23,\ 24,\ 25,\ 32 \; \}, \quad T = \{ \; 32 \; \}.
$$

Let us take an Ap\'ery ordering $\prec_1$, taking the lex ordering $y_1<y_2<y_3$ in step (3). The Groebner basis is
\begin{eqnarray*}
\cB_1 &=& \Big\{  x^7 - y_1, \ x^4 y_3 - y_4,\ x^2 y_3^2 - y_1 y_4,\ x y_2 - y_3,\ x y_1 - y_2,\ x y_4 - y_1^2, \\
&& \quad \quad y_2^5 - y_1^2 y_4^2,\ y_2^3 y_3 - y_1 y_4^2,\ y_1 y_2^3 - y_3^2 y_4,\ y_3^3 - y_1^2 y_4,\ y_2 y_3^2 - y_4^2, \\
&& \quad \quad y_1^2 y_2 - y_3 y_4,\ y_1^3 - y_2 y_4,\ y_1 y_3 - y_2^2 \Big\}.
\end{eqnarray*}

Let us compute $\Delta_{\prec_1} (S,13)$. In order to do that let us draw in $\{x=y_4=0\}$ the integer corresponding by $\cG^{-1}$ to each element $(0,\alpha_1,\alpha_2,\alpha_3,0)$ (that is, $7\alpha_1 + 8 \alpha_2 + 9 \alpha_3$). The shadowed regions in the picture below are precisely $\overline{E(I_S)} \cap \{ x=y_4=0 \}$ (the horizontal axis is $y_1$, the vertical is $y_2$). 

At $y_3=0$

\vspace{.5cm}

\includegraphics[scale=0.65]{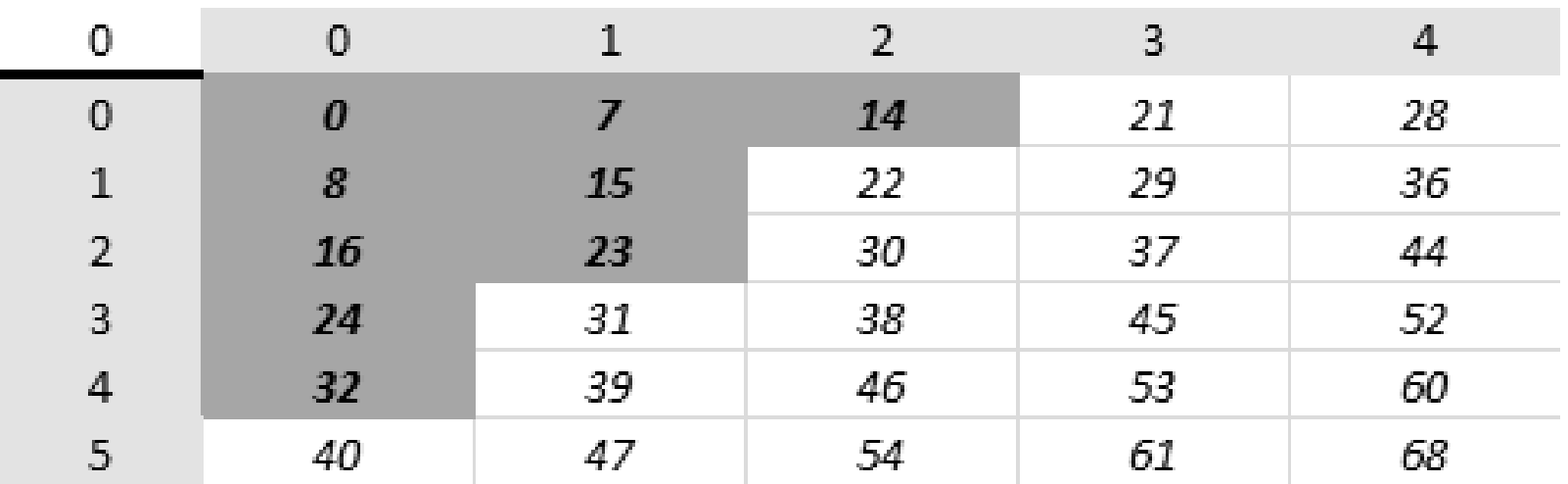}

\vspace{.5cm}

On the other hand, at $y_3=1$

\vspace{.5cm}

\includegraphics[scale=0.65]{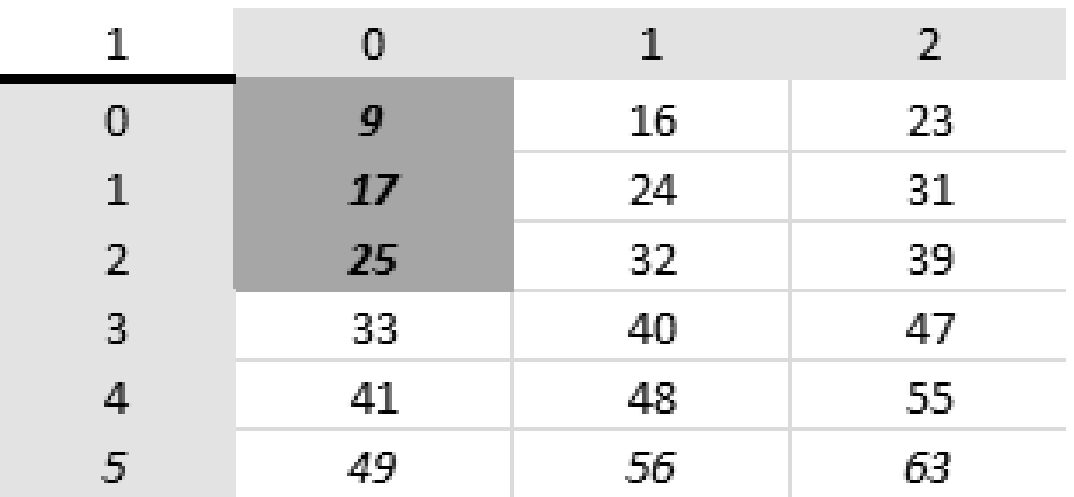}

\vspace{.5cm}

Finally, at $y_3=3$

\vspace{.5cm}

\includegraphics[scale=0.65]{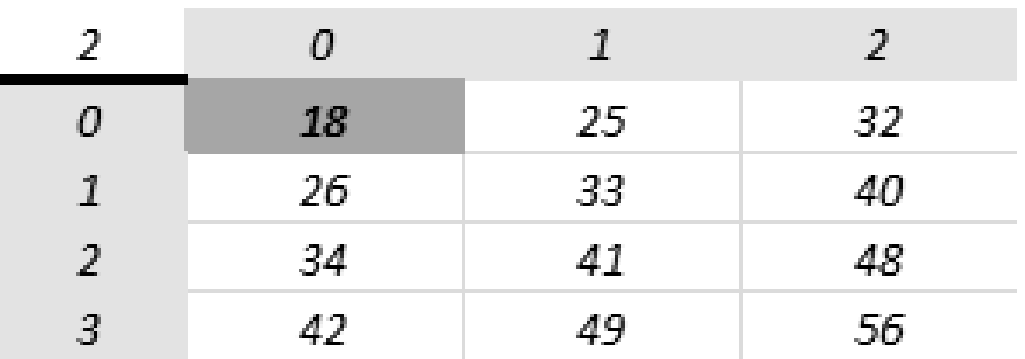}

\vspace{.5cm}

That is, as expected,
$$
\Delta_{\prec_3} (S,13) = \{ \; 0,\ 7,\ 8,\ 9,\ 14,\ 15,\ 16,\ 17,\ 18,\ 23,\ 24,\ 25,\ 32 \; \}.
$$

Mind that, although the set $\Delta_{\prec_j} (S,a_j)$ does not depend on the chosen ordering for $\{y_1,...,y_k\} \setminus \{y_j\}$, the actual arranging of the integers inside $\Delta_{\prec_j} (S,a_j)$ might well depend. In the same example, consider now in the step (3) the ordering $\prec_2$, taking the lex ordering $y_2<y_3<y_1$ in step (3). The Groebner basis is
\begin{eqnarray*}
\cB_2 &=& \Big\{  x^7 - y_1,\ x^4 y_3 - y_4,\ x^2 y_3^2 - y_1 y_4,\ x y_2 - y_3,\ x y_1 - y_2,\ x y_4 - y_1^2, \\
&& \quad \quad y_3^3 - y_1^2 y_4,\ y_2y_3^2 - y_4^2,\ y_1^2 y_2 - y_3 y_4,\ y_1^3 - y_2 y_4,\ y_2^2 - y_1 y_3 \Big\}.
\end{eqnarray*}

Then the picture of $\Delta_{\prec_3}(S,13)$ is, at $x=y_4=y_3=0$,

\vspace{.5cm}

\includegraphics[scale=0.65]{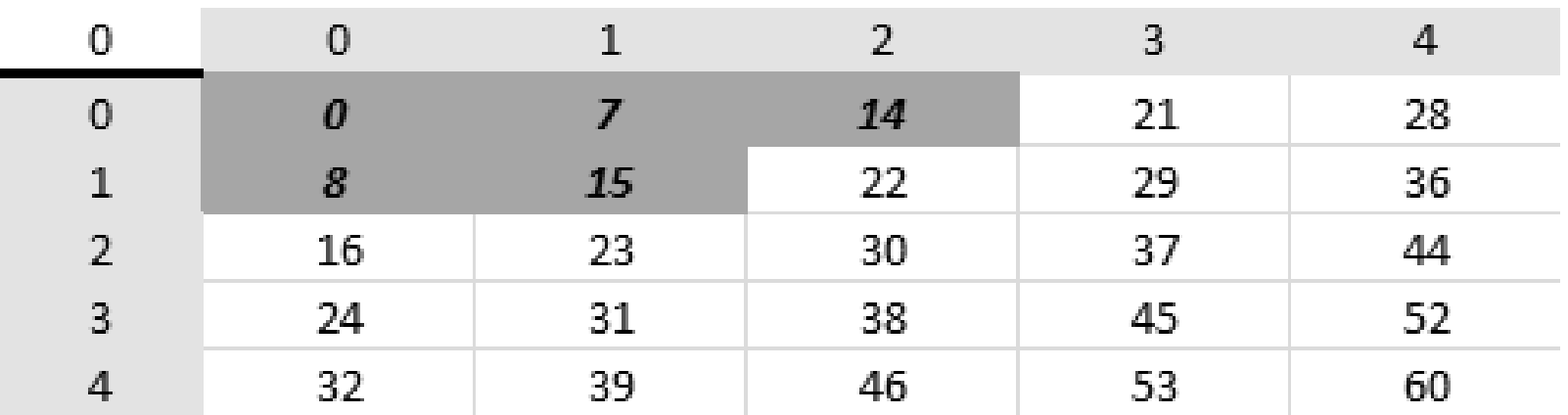}

\vspace{.5cm}

At $x=y_4=0, \; y_3=1$,

\vspace{.5cm}

\includegraphics[scale=0.65]{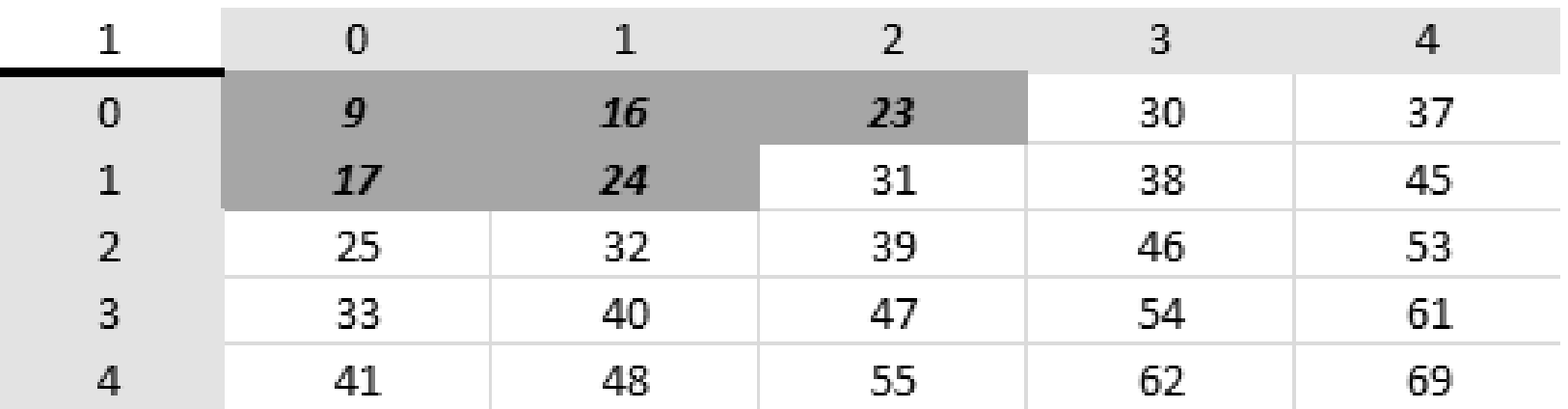}

\vspace{.5cm}

And at last, at $x=y_4=0, \; y_3=2$,

\vspace{.5cm}

\includegraphics[scale=0.65]{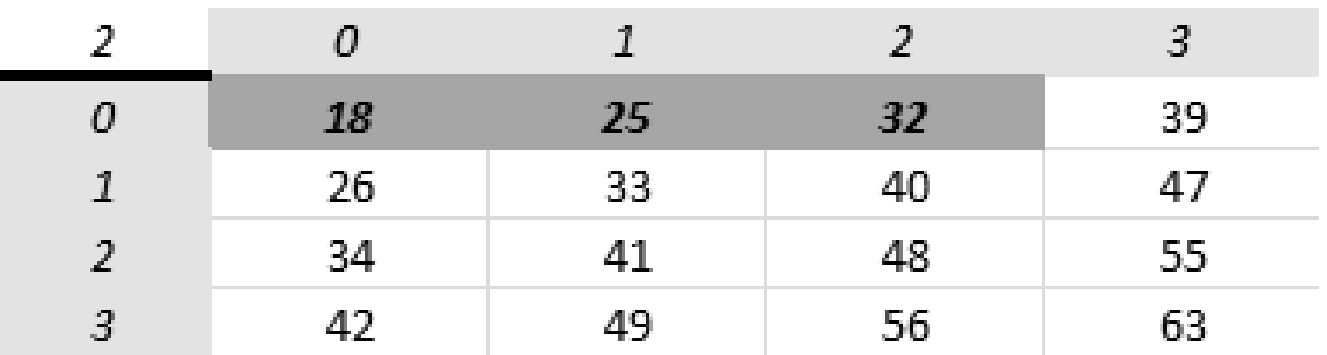}

\vspace{.5cm}
\end{example}

\section{Computation of the Ap\'ery set: the affine case}

An affine monoid is a finitely generated monoid that is isomorphic to a submonoid of $\mathbb{Z}^d, d \geq 0$. For the sake of simplicity, we will assume that all our affine monoids are submonoids of  $\mathbb{Z}^d$ for some $d \geq 0$.

Let $S = \langle \mathbf{a}_1, \ldots, \mathbf{a}_k \rangle$ be an affine monoid with $\mathbf{a}_i \in \mathbb{Z}^d \setminus \{0\}$. We will say that $S$ is pointed if $S \cap (-S) = \{0\}$, that is to say, if $\textbf{0}$ is the only invertible element of $S.$ Equivalently, if the rational cone 
$$
pos(S) := \{ \lambda_1 \mathbf{a}_1 + \ldots + \lambda_k \mathbf{a}_k \mid \lambda_i \in \mathbb{Q}_{\geq 0} \}
$$ 
is pointed. 

Observe that if $d=1$, then pointed affine monoids are nothing but numerical monoids. In this case, the corresponding rational cones are all equal to $\mathbb{Q}_{\geq 0}$.

Let $\Lambda \subseteq \{\mathbf{a}_1, \ldots, \mathbf{a}_k\}$ be such that $pos(S) = pos(\Lambda)$. The Ap\'ery set of $S$ with respect to $\Lambda$ is defined as follows (see, e.g., \cite[Definition 1.1]{Pilar}) 
$$
Ap(S,\Lambda) = \{ \mathbf{a} \in S \mid \mathbf{a} - \mathbf{b} \not\in S, \forall \mathbf{b} \in \Lambda\}.
$$ 
For $d = 1$, one has that $\mathrm{pos}(S) = \mathrm{pos}(s)$, for all $s \in S \setminus \{0\}$. Thus, by Lemma \ref{lema0}, the above definition of Ap\'ery set is a generalization of the one given for numerical semigroups.

Let us fix $\Lambda \subseteq \{\mathbf{a}_1, \ldots, \mathbf{a}_k\}$ such that $pos(S) = pos(\Lambda)$. Without loss of generality, we may suppose that $\Lambda = \{\mathbf{a}_{k-n+1}, \ldots, \mathbf{a}_k\},\ n \leq k$.

Consider the polynomial ring $\mathbb{Q}[x_1, \ldots, x_d,y_1, \ldots, y_k]$ and let $\prec_\Lambda$ be a block ordering on $A$ such that $\prec_\Lambda$ is an arbitrary monomial ordering $\prec_x$ for $x$ and $\prec_\Lambda$ is a $S$--graded reverse lexicographical ordering $\prec_{y}$ for $y$ such that $y_j \prec_y y_i$, for every $j \in \{1, \ldots, k-n\}$ and $j \in \{k-n+1, \ldots, k\}$.

\begin{example}
If $\prec_x$ is the lexicographical ordering with $x_d \prec_x \ldots \prec_x x_1$ and $\prec_y$ is the $S$--graded reverse lexicographical ordering with $y_k \prec_y \cdots \prec_y y_1$, then $\prec_\Lambda$ is given by 
$$
\left(\begin{array}{c|cccc} 
I_d & 0 & 0 & \cdots & 0 \\ 
\cline{1-5} 0 & a_1 & a_2 & \cdots & a_k \\ 
0 & 0 & 0 & \cdots & -1 \\ 
\vdots & \vdots & \vdots & & \vdots \\ 
0 & 0 & -1 & \cdots & 0  \end{array}\right)
$$
\end{example}

Let $I_S$ be the kernel of the ring homomorphism 
\begin{eqnarray*}
\widetilde\phi : \mathbb{Q}[x_1, \ldots, x_d,y_1, \ldots, y_k] & \longrightarrow & \mathbb{Q}[x_1, \ldots, x_d] \\ 
y_j & \longmapsto & \mathbf{x}^{\mathbf{a}_j} := x_1^{a_{1j}} \cdots x_d^{a_{dj}} \\ 
x_i & \longmapsto & x_i
\end{eqnarray*} 

Let $\mathcal{B}_\Lambda$ be the reduced Groebner basis of $I_S$ with respect to $\prec_\Lambda$ and set $N_\Lambda$ for the normal form operator with respect to this basis. For the sake of convenience, we will write $\{z_1, \ldots, z_n\}$ instead of $\{y_{k-n+1},\ldots, y_k\}$ in what follows.

Now, define the following set 
$$
\mathcal{Q}_{\prec_\Lambda}(S) = \Big\{ \mathbf{a} \in \mathbb{Z}^d_{\geq 0} \mid \exp\big(N_\Lambda(\mathbf{x}^\mathbf{a}) \big) \in \{x_1= \ldots x_d = z_1 = \ldots = z_n = 0\} \cap \overline{E(I_S)} \Big\}
$$

Notice that there are finitely many elements in $\mathcal{Q}_{\prec_\Lambda}(S)$. Indeed, since $pos(S) = pos(\Lambda)$, then, for each $a_j,\ j = 1, \ldots, k-n$, there exist $u_j \in \mathbb{Z}_{\geq 0}$ and $\mathbf{v}_j = (v_{1j}, \ldots, v_{nj}) \in \mathbb{Z}^n_{\geq 0}$ such that $u_j \mathbf{a}_j = \sum_{i=1}^n v_{ij} \mathbf{a}_{k-n+i}$. Therefore, 
$$
y_j^{u_j} - N_\Lambda(\mathbf{z}^{\mathbf{v}_j})  \in \mathcal{B}_\Lambda.
$$

The next theorem generalizes Theorem \ref{Apery} for affine monoids.

\begin{theorem}
$\mathcal{Q}_{\prec_\Lambda}(S) = Ap(S,\Lambda)$.
\end{theorem}

\begin{proof}
Let $\mathbf{a} \in \mathcal{Q}_{\prec_\Lambda}(S)$ and consider the fiber of  
\begin{eqnarray*}
\pi : \mathbb{Z}^{d+k}_{\geq 0} & \longrightarrow & S \\ 
(u_1, \ldots, u_{d+k}) & \longmapsto & \sum_{i={d+1}}^{d+k} u_i \mathbf{a}_k
\end{eqnarray*}
over $\mathbf{a}$, that will be denoted by $\pi^{-1}(\mathbf{a})$. 

By hypothesis, there exists 
$$
\mathbf{u} = (0, \ldots,0,u_{d+1},\ldots, u_{d+k-n}, 0, \ldots, 0) \in \pi^{-1}(\mathbf{a}) \cap \overline{E(I_S)}.
$$
If $\mathbf{a} - \mathbf{b} \in S,$ for some $\mathbf{b} \in \Lambda,$ there exist $\mathbf{v} = (0, \ldots, 0, v_{d+1}, \ldots, v_{d+k-n}, 0, \ldots, 0)$ and $\mathbf{w} = (0, \ldots, 0, w_1, \ldots, w_n) \in \mathbb{Z}^{d+k}_{\geq 0}$ with $w_i \neq 0$, for some $i$, and $\mathbf{v} + \mathbf{w} \in  \pi^{-1}(\mathbf{a})$. 

That is to say, there exists a nonzero 
$$
f = \mathbf{y}^\mathbf{u} - \mathbf{y}^\mathbf{v} \mathbf{z}^\mathbf{w} \in I_S.
$$

Finally, since, by the definition of $\prec_\Lambda$, $\mathbf{v} + \mathbf{w} \prec_\Lambda \mathbf{u}$, we obtain that $\mathbf{u} \in E(I_S)$, a contradiction.

Conversely, if $\mathbf{a} \in Ap(S,\Lambda)$ and $\mathbf{u} \in \pi^{-1}(\mathbf{a})$, then $u_i = 0,$ for every $i \in \{{d+k-n+1}, \ldots, {d+k}\}$. In particular, 
$$
\exp(N_\Lambda(\mathbf{x}^\mathbf{a})) \in \{x_1= \ldots x_d = z_1 = \ldots = z_n = 0\}
$$ 
and we are done.
\end{proof}

\section{Pseudo--Frobenius numbers and the type set}

Let us fix a numerical monoid $S$ and let us review some classical definitions. First, we define a partial ordering $\leq_S$ in $S$ as follows: 
$$
x \leq_S y \; \Longleftrightarrow \; y-x \in S.
$$

We also say $x \in \Z$ is a  pseudo--Frobenius number for $S$ if:
\begin{itemize}
\item $x \notin S$.
\item $x+s \in S$ for all $s \in S \setminus \{0\}$. 
\end{itemize}

The set of pseudo--Frobenius numbers for $S$ will be noted $PF(S)$ and its cardinal, which we will call the type of $S$ will be written as $t(S)$.

It is clear that $f(S)$ is the maximum of $PF(S)$ (with respect to the usual ordering in $\Z$). A very special familiy of monoids are precisely those where $PF(S) = \{ f(S) \}$ (equivalently, where $t(S)=1$). These monoids are called symmetric. 

The relationship between $<_S$ and $PF(S)$ comes from the next result.

\begin{proposition}
An integer $g \in \Z$ belongs to $PF(S)$ if and only if for any $n \in S$, $g + n$ is a maximal element in $Ap(S,n)$ with respect to the ordering $\leq_S$. 
\end{proposition}

\begin{corollary}
$t(S) = \# \left\{ \max_{\leq_S} Ap(S,n) \right\}.$
\end{corollary}

In \cite{NW} Nijenhuis and Wilf studied when the property 
$$
g(S) = \frac{1}{2} \left( f(S)+1 \right)
$$
holds for a given numerical semigroup $S$ (the same goes for monoids, of course). As we mentioned, the case $k=2$ was completely known from Sylvester \cite{Sylvester}. Looking for conditions that are equivalent to this property, they proved the following result (essentially proved independently by Kunz \cite{Kunz}).

\begin{theorem}
Let $S = \langle a_1,...,a_k \rangle$ and let us consider the set
$$
T(S) = \left\{ m \in Ap(S,a_k) \; | \; m+a_i \notin Ap(S,a_k), \; \forall i=1,...,k \right\}.
$$

Then 
$$
g(S) = \frac{1}{2} \left( f(S)+1 \right) \; \Longleftrightarrow \; \sharp T(S) = 1.
$$
\end{theorem}

The set $T(S)$ will be called the type set of $S$, and the condition $\sharp T(S) =1$ will be called (after \cite{NW}) Gorenstein condition. 

There is also a tight relationship between $T(S)$ and $PF(S)$ \cite{NW}.

\begin{proposition}
$PF(S) = \left\{ m-a_k \; | \; m \in T(S) \right\}$.
\end{proposition}

\begin{corollary}
Under the previous assumptions:
\begin{itemize}
\item $\# \ T(S) = t(S)$.
\item A monoid verifies the Gorenstein condition if and only if it is symmetric.
\item $T(S)= \big\{ \max_{\leq_S} Ap(S,a_k) \ \big\}$.
\end{itemize}
\end{corollary}

In the following section we will give a different method for computing the type set using the algorithm we developed above.

\section{Computation of the type set}

\begin{example} 
Let us take again our monoid from a previous example $S = \langle 3,7,11 \rangle$. Remember we already computed
$$
Ap(S,11)= \{ 0, \ 3, \ 6, \ 7, \ 9, \ 10, \ 12, \ 13, \ 15, \ 16, \ 19 \}.
$$

If we want to compute the type set $T(S)$ and the set $PF(S)$ we can use the partial ordering $\leq_S$ as indicated above. In fact, the ordering induces the following diagram in $Ap(S,11)$ (an arrow $x \longrightarrow y$ indicates $x \leq_S y$):

\vspace{.3cm}

\begin{center}
\quad \quad \xymatrix{
 & 7 \ar[r] & 10 \ar[r] & 13 \ar[r] & 16 \ar[r] & 19\\
 0 \ar[ur] \ar[r] & 3 \ar[ur]\ar[r] & 6 \ar[ur] \ar[r] & 9 \ar[ur]\ar[r] & 12 \ar[ur]\ar[r] & 15}
\end{center}

\vspace{.3cm}

Therefore we have $T(S) = \{15, \ 19\}$, $PF(S) = \{4,\ 8 \}$.
\end{example}

In this section we will give a different, systematic way of computing $T(S)$, using our previous algorithm to compute $Ap(S,s)$.

Under the assumptions of the previous section, we compute $Ap(S,a_k)$ using an Ap\'ery monomial ordering $\prec$ whose normal form will be denoted by $N_k(\cdot)$. Let $N \in Ap(S,a_k)$, and therefore let us write
$$
\exp \left( N_k \left( x^N \right) \right) = (0,\gamma_1,...,\gamma_{k-1},0).
$$

We will say $N$ is an extremal element of $Ap(S,a_k)$ for $\prec$ if, for all $i=1,...,k$ we have
$$
(0,\gamma_1,...,\gamma_i+1,...,\gamma_{k-1},0) \notin \overline{E(I_S)} \cap \{ \ x=y_k=0 \ \}.
$$

Mind that the set of extremal elements {\em does} depend on the ordering chosen $\prec$. That is why we will denote
$$
\partial_{\prec} (S,a_k) = \left\{ \mbox{ Extremal elements of $Ap(S,a_k)$ for $\prec$ } \right\}.
$$

\begin{lemma}
$T(S) \subset \partial_{\prec} (S,a_k)$ for any Ap\'ery ordering $\prec$ with respect to $a_k$.
\end{lemma}

\begin{proof}
Assume it is not so, for some $N \in T(S)$ with
$$
\exp \left( N_k \left( x^N \right) \right) = (0,\gamma_1,...,\gamma_{k-1},0). 
$$
Then there exists $i \in \{1,...,k-1\}$ such that
$$
(0,\gamma_1,...,\gamma_i+1,...,\gamma_{k-1},0)  \in \overline{E(I_S)} \cap \{ \ x=y_k=0 \ \}.
$$

This implies, from Theorem \ref{Apery}, and using $\cG^{-1}$ from Theorem \ref{GC}, that $N+a_i \in Ap(S,a_k)$, contradicting the fact that $N \in T(S)$.
\end{proof}

\begin{theorem}\label{GS}
If $\cO$ is the set of Ap\'ery orderings with respect to $a_k$, then
$$
\bigcap_{\prec \in \cO} \partial_\prec (S,a_k) = T(S).
$$
\end{theorem}

\begin{proof}
By the previous lemma, we only need to proof the following: If $N \notin T(S)$, then there exists an Ap\'ery ordering $\prec \in \cO$ such that $N \notin \partial_\prec (S,a_k)$.

If $N \notin T(S)$, there must exist a generator $a_i$ for $i=1,...,k-1$ such that $N+a_i \in Ap(S,a_k)$. Let us take $\prec \in \cO$ any monomial ordering such that, in step (3), we choose the reverse lexicographic ordering with respect to 
$$
\left\{ y_1,...,y_{i-1}, y_{i+1},..., y_{k-1},y_i \right\}.
$$

That is, given two $(k+1)$--uples, the smallest is the one with bigger $(i+1)$--th coefficient (corresponding with $y_i$) and so on. Mind that this is a monomial ordering (in particular, it is a well ordering) only beacuse of the grading stablished in step (2) of the definition of Ap\'ery orderings.

Let us write
$$
\exp \left( N_\prec \left( x^{N} \right) \right) = (0,\beta_1,...,\beta_i,...,\beta_{k-1},0).
$$

Then it must hold that, 
$$
\exp \left( N_\prec \left( x^{N+a_i} \right) \right) = (0,\beta_1,..., \beta_i+1, ..., \beta_{k-1},0), 
$$

This comes straightforwardly, if we had otherwise
$$
\exp \left( N_\prec \left( x^{N+a_i} \right) \right) = (0,\alpha_1,..., \alpha_i, ..., \alpha_{k-1},0), 
$$
then, from our ordering, we must have either $\alpha_i > \beta_i+1$, or $\alpha_i = \beta_i+1$ and we break the tie in another coordinate of the $(k+1)$--uple. 

This last case cannot happen, because the rest of the coordinates must be the coefficients of the smallest representation (with respect to our chosen ordering) of $N-\beta_ia_i$.

The first option is also impossible, because if $\alpha_i > \beta_i+1$, then it is clear that
$$
(0,\alpha_1,..., \alpha_i-1, ..., \alpha_{k-1},0) \in \Z_{\geq 0}^k
$$
is a $(k+1)$--uple which represents $N$ and it is smaller than $(0,\beta_1,...,\beta_i,...,\beta_{k-1},0)$ with respect to our ordering. Therefore $N \notin \partial_\prec (S,a_k)$.
\end{proof}

Despite our result displays a infinite intersection it is clear from the proof that we need, at most, $k-1$ monomial orderings (and the subsequent computations of the sets $\partial_\prec (S,a_k)$) for the full computation of $T(S)$.

\begin{example} 
Let us consider again the monoid $S = \langle 7,\ 8,\ 9,\ 13 \rangle$, which verifies
$$
G(S)= \{ \; 1, \ 2, \ 3, \ 4, \ 5, \ 6, \ 10, \ 11, \ 12, \ 19 \; \}
$$
and 
$$
Ap(S,13) = \{ \; 0, \ 7, \ 8, \ 9, \ 14, \ 15, \ 16, \ 17, \ 18, \ 23, \ 24, \ 25, \ 32 \},
$$
$$
T(S) = \{ \; 32 \; \},
$$
that is, $S$ is symmetric.

From Theorem \ref{GS} we need to consider (at most) $3$ monomial orderings. Let us write $\prec_1$ the Ap\'ery ordering which takes, in step (3), the reverse lexicographic ordering with respect to $\{ y_1, \ y_2, \ y_3 \}$. The Grobner basis is precisely the basis $\cB_2$ in Example 3.2. From the arrangement of the set $Ap(S,13)$ according to this basis, we see that the extremal elements are precisely:
$$
\partial_{\prec_1} (S,13) = \{ \; 24, \ 32 \}.
$$

Note that $24 \notin T(S)$ because $24 + a_2 = 24+8 \in Ap(S,13)$. The same results are obtained if we use the ordering $\prec_3$, which uses the reverse lexicographic ordering with respect to $\{ y_2, \ y_3, \ y_1 \}$.

However, if we consider (following the strategy of the proof) the ordering $\prec_2$, taking in step (3) the reverse lexicographic ordering with respect to $\{ y_1, \ y_3, \ y_2 \}$, we get the Grobner basis $\cB_1$ from Example 3.2.

In this case we have
$$
\partial_{\prec_2} (S,13) = \{ \; 14,\ 18,\ 23,\ 25,\ 32 \; \},
$$
so we have four stowaways, namely $\{14,\ 18,\ 23,\ 25 \}$ which verify
$$
14 + a_3 = 14+9 \in Ap(S,13), \quad 18 + a_1 = 18+7 \in Ap(S,13),
$$
$$
23 + a_3 = 23+9 \in Ap(S,13), \quad 25 + a_1 = 25+7 \in Ap(S,13).
$$

As expected, $T(S) = \{\ 32 \ \} = \partial_{\prec_1} (S,13) \cap \partial_{\prec_2} (S,13)$.
\end{example}

\section{Another interpretation of pseudo--Frobenius numbers}

We will end with a word on an important feature of the pseudo--Frobenius elements (and hence the set $T(S)$) in a more highbrow context. Let 
$$
S = \langle a_1, \ldots, a_k \rangle \subset \mathbb{Z}_{\geq 0}
$$ 
be a numerical monoid and let $\mathbb{Q}[S]$ be the monoid algebra associated to $S$, that is to say, 
$$
\mathbb{Q}[S] := \bigoplus_{s \in S} \mathbb{Q}\, \chi^s,\ \text{with}\ \chi^s \cdot \chi^{s'} := \chi^{s+s'}.
$$

Obviously, the monoid algebra $\mathbb{Q}[S]$ is an $S-$graded ring. We also consider the polynomial ring $A := \mathbb{Q}[y_1, \ldots, y_k]$ as an $S$--graded ring, by assigning degree $a_i$ to $y_i$. Thus, the ring homomorphism 
\begin{eqnarray*}
\varphi_0 : A & \longrightarrow & \mathbb{Q}[S] \\ 
y_i & \longmapsto & \chi^{a_i}
\end{eqnarray*}
is clearly $S$--graded too. 

This allows to construct an $S$--graded minimal free resolution of $\mathbb{Q}[S]$ as $A$--mo\-du\-le that will be finite of length $\dim(A) - \mathrm{depth}(\mathbb{Q}[S]) = k - 1$, by the Auslander-Buchbaum theorem. More explicitly, the resolution has the form 
$$
0 \to \bigoplus_{j=1}^{l_{k-1}} A^{r_{k-1\, j}}(-b_{k-1\, j}) \longrightarrow \ldots \longrightarrow \bigoplus_{j=1}^{l_1} A^{r_{1\, j}}(-b_{1\, j}) \longrightarrow A \stackrel{\varphi_0}{\longrightarrow} \mathbb{Q}[S] \to 0,
$$ 
where the integer $\beta_i := \sum_{j=1}^{l_j} r_{ij}$ is the rank of the $i$--th syzygy module of $\mathbb{Q}[S]$ and $b_{ij} \neq b_{ij'},\ j \neq j'.$

The integers $b_{ij}$ do not depend on the resolution and they can be combinatorially characterized as follows.  

Let $a \in S$. Consider the abstract simplicial complex $\Delta_a$ consisting of all subsets $F$ of $\{1, \ldots, k\}$ such that $a - \sum_{i \in F} a_i \in S,$ and let $\widetilde{H}_i(\Delta_a)$ be the $i$--th simplicial homology vector space of $\Delta_a$ with values in $\mathbb{Q}$.

\begin{theorem}
$\widetilde{H}_{i-1}(\Delta_a) \neq 0$ if and only if $a = b_{ij}$ for some $j$. Moreover, in this case, $\dim(\widetilde{H}_{i-1}(\Delta_a)) = r_{ij}$
\end{theorem}

\begin{proof}
See \cite[Theorem 2.1]{CM}.
\end{proof}

\begin{corollary}
$PF(S) = \big\{ b - \sum_{i=1}^k a_i \mid b \in \{b_{k-1\, 1}, \ldots,  b_{k-1\, l_{k-1}} \} \big\}$.
\end{corollary}

\begin{proof}
It suffices to observe that $a \in PF(S)$ if and only if $\Delta_{a + \sum_{i=1}^k a_i}$ has the reduced homology of a $(k-2)$-sphere.
\end{proof}

\section{Final remarks}

The computation of the Ap\'ery set and the type set may be seen as a first step in the understanding of more complicated structures inside a numerical monoid. In particular, different partial orderings from $\leq_S$ can be considered in order to gather more interesting data (see for instance \cite{DMS}). We hope this work may serve as a first approach to tackle this further--reaching problems.

The authors wish to thank P. Garc\'{\i}a--S\'anchez, for pointing us the example which eventually led to Theorem \ref{GS}, to J.L. Ram\'{i}rez--Alfons\'{\i}n and M. D'Anna for their help and advice and also to S. Robbins for his enlightening conversations during his stay in Seville in December 2013.

\end{document}